\documentclass[10pt]{amsart}
\usepackage[psamsfonts]{amssymb}
\usepackage{amsthm,amscd}
\usepackage{type1cm}
\newtheorem{thm}{Theorem}[section]

\theoremstyle{remark}

\theoremstyle{definition}

\allowdisplaybreaks[4]

\title[A characterization of endo-commutativity]{A characterization of endo-commutativity of 3-dimensional curled algebras}

\author[S.-E. Takahasi]{Sin-Ei Takahasi}
\author[K. Shirayanagi]{Kiyoshi Shirayanagi}
\address[S.-E. Takahasi]{Laboratory of Mathematics and Games\\ Innaityo 668-2\\ Funabashi\\ Chiba 273-0025\\ Japan}
\address[K. Shirayanagi]{Department of Information Science\\ Toho University\\ Miyama 2-2-1\\ Funabashi\\ Chiba 274-8510\\ Japan}
\email{sin\_ei1@yahoo.co.jp}
\email[K. Shirayanagi (Corresponding author)]{kiyoshi.shirayanagi@is.sci.toho-u.ac.jp}

\makeatletter
\@namedef{subjclassname@2020}{\textup{2020} Mathematics Subject Classification}
\makeatother

\subjclass[2020]{Primary 17A30; Secondary 17D99, 13A99}
\keywords{Nonassociative algebras, curled algebras, endo-commutative algebras}

\date{\today}

\begin{document}

\begin{abstract}
  A curled algebra is a non-associative algebra in which $x$ and $x^2$ are linearly dependent for every element $x$.
  An algebra is called endo-commutative, if the square mapping from the algebra to itself preserves multiplication.
  In this paper, we provide a necessary and sufficient condition for a 3-dimensional curled algebra over an arbitrary field
  to be endo-commutative, expressed in terms of the properties of its underlying linear basis.
  
\end{abstract}

\maketitle

\section{Introduction}\label{sec:intro} 
Since the concept of zeropotent algebras was introduced as a generalization of exterior algebras or Lie algebras by the authors, along with Kobayashi and Tsukada,
it has been cited by many papers, including \cite{zeropotent 1}, \cite{zeropotent 2}, \cite{zeropotent 3}, \cite{AANS}, and \cite{FF}.

Subsequently, the concept of endo-commutativity was introduced as a further generalization of zeropotency by the authors and Tsukada a few years ago.
Following this, a complete classification of 2-dimensional endo-commutative curled algebras was successfully provided, as detailed in \cite{TST}.

In this paper, we provide a necessary and sufficient condition for a 3-dimensional curled algebra over an arbitrary field
to be endo-commutative, expressed in terms of the properties of its underlying linear basis. This result is seen as a step toward
solving the challenging problem of classifying 3-dimensional endo-commutative curled algebras.

\vspace{3mm}

\section{3-dimensional endo-commutative curled algebras}\label{sec:ec-curled} 
Let $K$ be an arbitrary field and let $\mathcal A$ be a (not necessarily associative) algebra over $K$.  An element $x\in \mathcal A$ is said to be {\it{curled}}, if the set $\{x,  x^2\}$ is linearly dependent over $K$.  The algebra $\mathcal A$ is said to be curled if every element of $\mathcal A$ is curled. 

Suppose $\mathcal A$ is a 3-dimensional curled algebra over $K$ with a linear basis $\{e, f, g\}$.  Since $e, f$ and $g$ are all curled,  there exist scalars $\varepsilon_e, \varepsilon_f, \varepsilon_g\in K$ such that
\[e^2=\varepsilon_ee, f^2=\varepsilon_ff, \mbox{ and } g^2=\varepsilon_gg.\]
By a suitable change of basis, we may assume without loss of generality that  $\varepsilon_e, \varepsilon_f, \varepsilon_g\in \{0, 1\}$.  Given a triple $(i, j, k)\in \{0, 1\}^3$, we say that $\mathcal A$ is a 3-dimensional curled algebra of
{\it{type {\rm{(}}i, j, k{\rm{)}}}} with respect to the basis $\{e, f, g\}$ if $(\varepsilon_e, \varepsilon_f, \varepsilon_g)=(i, j, k)$.
Thus, there are exactly eight possible types of 3-dimensional curled algebras.

An algebra $\mathcal A$ over $K$ is said to be {\it{endo-commutative}} if the square mapping from $\mathcal A$ to itself preserves multiplication; that is, for all $x, y\in \mathcal A$, we have \[x^2y^2=(xy)^2.\]

Let  $i,j, k\in \{0, 1\}$ and suppose that $\mathcal A$ is a 3-dimensional curled algebra of type $(i, j, k)$.   Thus,  $\mathcal A$ has a multiplication table 
\[
\begin{pmatrix} ie&ef&eg\\fe&jf&fg\\ge&gf&kg\end{pmatrix}.
\]
Define
\[
A=ef, B=eg, C=fe, D=fg, E=ge,\, \, {\rm{and}}\, \, F=gf.
\]

In the next section, we provide a necessary and sufficient condition for a 3-dimensional curled algebra of type $(i, j, k)$ to be endo-commutative,
in terms of the relations among the basis elements $\{e, f, g\}$ and
the parameters $A, B, C, D, E, F$.
\vspace{3mm}

\section{A characterization of endo-commutativity of curled algebras}\label{sec:characterization} 
Let $x, y\in \mathcal A$ be arbitrary elements, and write \[x=ae+bf+cg, y=ue+vf+wg,\] where $a, b, c, u, v, w\in K$.  Then
\[
xy=auie+avA+awB+buC+bvjf+bwD+cuE+cvF+cwkg.
\]
We define
\begin{align*}
[(xy)^2]_1&:=a^2u^2ie+a^2v^2A^2+a^2w^2B^2+b^2u^2C^2+b^2v^2jf+b^2w^2D^2\\
&\quad\quad+c^2u^2E^2+c^2v^2F^2+c^2w^2kg,
\end{align*}
\begin{align*}
[(xy)^2]_2&:=a^2uvi(eA+Ae)+a^2uwi(eB+Be)+abu^2i(eC+Ce)\\
&\quad\quad+abuvij(A+C)+abuwi(eD+De)+acu^2i(eE+Ee)\\
&\quad\quad+acuvi(eF+Fe)+acuwik(B+E),
\end{align*}
\begin{align*}
[(xy)^2]_3&:=a^2vw(AB+BA)+abuv(AC+CA)+abv^2j(Af+fA)\\
&\quad\quad +abvw(AD+DA)+acuv(AE+EA)+acv^2(AF+FA)\\
&\quad\quad +acvwk(Ag+gA),
\end{align*}
\begin{align*}
[(xy)^2]_4&:=abuw(BC+CB)+abvwj(Bf+fB)+abw^2(BD+DB)\\
&\quad\quad+acuw(BE+EB)+acvw(BF+FB)+acw^2k(Bg+gB),
\end{align*}
\begin{align*}
[(xy)^2]_5&:=b^2uvj(Cf+fC)+b^2uw(CD+DC)+bcu^2(CE+EC)\\
&\quad\quad +bcuv(CF+FC)+bcuwk(Cg+gC),
\end{align*}
\begin{align*}
[(xy)^2]_6&:=b^2vwj(fD+Df)+bcuvj(fE+Ef)+bcv^2j(fF+Ff)\\
&\quad\quad+bcvwjk(D+F),
\end{align*}
\[
[(xy)^2]_7:=bcuw(DE+ED)+bcvw(DF+FD)+bcw^2k(Dg+gD),
\]
\[
[(xy)^2]_8:=c^2uv(EF+FE)+c^2uwk(Eg+gE)
\]
and
\[
[(xy)^2]_9:=c^2vwk(Fg+gF).
\]
By simple calculations,  we obtain
\begin{equation} 
(xy)^2=\sum_{\ell=1}^9[(xy)^2]_\ell.
\end{equation}
Also, we have
\[
\left\{\begin{array}{@{\,}lll} 
x^2=a^2ie+b^2jf+c^2kg+ab(A+C)+ac(B+E)+bc(D+F),\\
y^2=u^2ie+v^2jf+w^2kg+uv(A+C)+uw(B+E)+vw(D+F).
\end{array} \right.
\]

We define
 \begin{align*}
[x^2y^2]_1&:=a^2u^2ie+a^2v^2ijA+a^2w^2ikB+a^2uvi(eA+eC)\\
&\quad\quad+a^2uwi(eB+eE)+a^2vwi(eD+eF),
\end{align*}
 \begin{align*}
[x^2y^2]_2&:=b^2u^2ij C+b^2v^2jf+b^2w^2jkD+b^2uvj(fA+fC)\\
&\quad\quad +b^2uwj(fB+fE)+b^2vwj(fD+fF),
\end{align*}
 \begin{align*}
[x^2y^2]_3&:=c^2u^2ikE+c^2v^2jkF+c^2w^2kg+c^2uvk(gA+gC)\\
&\quad\quad+c^2uwk(gB+gE)+c^2vwk(gD+gF),
\end{align*}
 \begin{align*}
[x^2y^2]_4&:=abu^2i(Ae+Ce)+abv^2j(Af+Cf)+abw^2k(Ag+Cg)\\
&\quad\quad+abuv(A+C)^2+abuw(A+C)(B+E)+abvw(A+C)(D+F),
\end{align*}
 \begin{align*}
[x^2y^2]_5&:=acu^2i(Be+Ee)+acv^2j(Bf+Ef)+acw^2k(Bg+Eg)\\
&\quad\quad +acuv(B+E)(A+C)+acuw(B+E)^2+acvw(B+E)(D+F),
\end{align*}
and
 \begin{align*}
[x^2y^2]_6&:=bcu^2i(De+Fe)+bcv^2j(Df+Ff)+bcw^2k(Dg+Fg)\\
&\quad\quad+bcuv(D+F)(A+C)+bcuw(D+F)(B+E)+bcvw(D+F)^2.
\end{align*}
By simple calculations,  we obtain
\begin{equation} 
x^2y^2=\sum_{\ell=1}^6[x^2y^2]_\ell.
\end{equation}

By examining equations (1) and (2), we list the coefficients of each term
in the expressions $(xy)$, $x^2y^2$, and $(xy)^2-x^2y^2$:
\vspace{3mm}

$\quad\quad(xy)^2$\quad\quad\quad  $x^2y^2$  $\quad\quad\quad\quad\quad(xy)^2-x^2y^2$
\vspace{2mm} 
      
$e :\quad a^2u^2i, \quad\quad\, \,   a^2u^2i$;$\quad\quad\quad\quad\quad\quad\quad  0$

$ f : \quad b^2v^2j, \quad\quad\, \,   b^2v^2j$;$\quad\quad\quad\quad\quad\quad\quad  0$

$g : \quad c^2w^2k \quad\quad\, \,  c^2w^2k$;$\quad\quad\quad\quad\quad\quad\quad  0$

$A : \quad abuvij, \quad\quad  a^2v^2ij$;$\quad\quad\quad\quad av(bu-av)ij$

$B : \quad acuwik, \quad\quad  a^2w^2ik$;$\quad\quad\quad\quad  aw(cu-aw)ik $ 

$C : \quad abuvij, \quad\quad b^2u^2ij$;$\quad\quad\quad\quad  bu(av-bu)ij $  

$D : \quad bcvwjk, \quad\quad b^2w^2jk$;$\quad\quad\quad\quad  bw(cv-bw)jk $ 

$E : \quad acuwik, \quad\quad c^2u^2ik$;$\quad\quad\quad\quad  cu(aw-cu)ik $  

$F : \quad bcvwjk, \quad\quad c^2v^2jk$;$\quad\quad\quad\quad  cv(bw-cv)jk $ 

$eA : \quad a^2uvi, \quad\quad  a^2uvi$;$\quad\quad\quad\quad\quad\quad\quad  0 $ 

$eB : \quad a^2uwi, \quad\quad a^2uwi,$;$\quad\quad\quad\quad\quad\quad\,   0 $  

$eC : \quad abu^2i, \quad\quad  a^2uvi$;$\quad\quad\quad\quad\quad  au(bu-av)i $ 

$eD : \quad abuwi, \quad\quad  a^2vwi$;$\quad\quad\quad\quad \quad  aw(bu-av)i $ 

$eE : \quad acu^2i, \quad\quad  a^2uwi$;$\quad\quad\quad\quad\quad  au(cu-aw)i $ 

$eF : \quad acuvi, \quad\quad  a^2vwi$;$\quad\quad\quad\quad\quad av(cu-aw)i  $ 

$fA : \quad abv^2j, \quad\quad  b^2uvj$;$\quad\quad\quad\quad\quad  bv(av-bu)j $ 

$fB : \quad abvwj, \quad\quad b^2uwj$;$\quad\quad\quad\quad\quad  bw(av-bu)j $ 

$fC : \quad b^2uvj, \quad\quad b^2uvj$; $\quad\quad\quad\quad\quad\quad\quad  0 $

$fD : \quad b^2vwj, \quad\quad b^2vwj$;$\quad\quad\quad\quad\quad\quad\quad  0 $ 

$fE : \quad bcuvj, \quad\quad b^2uwj$;$\quad\quad\quad\quad\quad  bu(cv-bw)j $ 

$fF : \quad bcv^2j,  \quad\quad b^2vwj$;$\quad\quad\quad\quad\quad  bv(cv-bw)j $

$gA : \quad acvwk, \quad\quad  c^2uvk$;$\quad\quad\quad\quad\quad  cv(aw-cu)k $

$gB : \quad acw^2k, \quad\quad  c^2uwk$;$\quad\quad\quad\quad\quad  cw(aw-cu)k $

$gC : \quad bcuwk, \quad\quad  c^2uvk$;$\quad\quad\quad\quad\quad  cu(bw-cv)k $ 

$gD : \quad bcw^2k, \quad\quad  c^2vwk$;$\quad\quad\quad\quad\quad  cw(bw-cv)k $

$gE : \quad c^2uwk, \quad\quad c^2uwk$;$\quad\quad\quad\quad\quad\quad\quad  0 $

$gF : \quad c^2vwk, \quad\quad  c^2vwk$;$\quad\quad\quad\quad\quad\quad\quad 0  $

$Ae : \quad a^2uvi, \quad\quad abu^2i$;$\quad\quad\quad\quad\quad  au(av-bu)i $ 

$Be : \quad a^2uwi, \quad\quad acu^2i$;$\quad\quad\quad\quad\quad  au(aw-cu)i $ 

$Ce : \quad abu^2i, \quad\quad abu^2i$;$\quad\quad\quad\quad\quad\quad\quad  0 $ 

$De : \quad abuwi, \quad\quad  bcu^2i$;$\quad\quad\quad\quad\quad  bu(aw-cu)i $ 

$Ee : \quad acu^2i, \quad\quad  acu^2i$;$\quad\quad\quad\quad\quad\quad\quad  0 $ 

$Fe : \quad acuvi, \quad\quad  bcu^2i$;$\quad\quad\quad\quad\quad  cu(av-bu)i $

$Af : \quad abv^2j, \quad\quad  abv^2j$;$\quad\quad\quad\quad\quad\quad\quad  0 $ 

$Bf : \quad abvwj, \quad\quad acv^2j$;$\quad\quad\quad\quad\quad  av(bw-cv)j $ 

$Cf :\quad b^2uvj, \quad\quad abv^2j$;$\quad\quad\quad\quad\quad  bv(bu-av)j $ 

$Df : \quad b^2vwj, \quad\quad  bcv^2j$;$\quad\quad\quad\quad\quad  bv(bw-cv)j $

$Ef : \quad bcuvj, \quad\quad acv^2j$;$\quad\quad\quad\quad\quad  cv(bu-av)j $ 

$Ff :\ \quad bcv^2j, \quad\quad  bcv^2j$;$\quad\quad\quad\quad\quad\quad\quad 0  $

$Ag : \quad acvwk, \quad\quad  abw^2k$;$\quad\quad\quad\quad\quad  aw(cv-bw)k $ 

$Bg : \quad acw^2k, \quad\quad  acw^2k$;$\quad\quad\quad\quad\quad\quad\quad  0 $

$Cg : \quad bcuwk, \quad\quad  abw^2k$;$\quad\quad\quad\quad\quad  bw(cu-aw)k $ 

$Dg : \quad bcw^2k, \quad\quad  bcw^2k$;$\quad\quad\quad\quad\quad\quad\quad  0 $ 
 
$Eg : \quad c^2uwk, \quad\quad  acw^2k$;$\quad\quad\quad\quad\quad  cw(cu-aw)k $ 
 
$Fg : \quad c^2vwk, \quad\quad  bcw^2k$;$\quad\quad\quad\quad\quad  cw(cv-bw)k $

$A^2 : \quad a^2v^2, \quad\quad  abuv$;$\quad\quad\quad\quad\quad\quad  av(av-bu) $

$B^2 : \quad a^2w^2, \quad\quad  acuw$;$\quad\quad\quad\quad\quad\quad  aw(aw-cu) $ 

$C^2 : \quad b^2u^2, \quad\quad  abuv$;$\quad\quad\quad\quad\quad\quad bu(bu-av)  $ 

$D^2 : \quad b^2w^2, \quad\quad  bcvw$;$\quad\quad\quad\quad\quad\quad  bw(bw-cv) $ 

$E^2 : \quad c^2u^2, \quad\quad  acuw$; $\quad\quad\quad\quad\quad\quad  cu(cu-aw) $

$F^2 : \quad c^2v^2, \quad\quad   bcvw$;$\quad\quad\quad\quad\quad\quad  cv(cv-bw) $

$AB : \quad a^2vw, \quad\quad  abuw$;$\quad\quad\quad\quad\quad  aw(av-bu) $ 

$AC : \quad abuv, \quad\quad  abuv$;$\quad\quad\quad\quad\quad\quad\quad  0 $

$AD : \quad abvw, \quad\quad  abvw$;$\quad\quad\quad\quad\quad\quad\quad  0 $ 

$AE : \quad acuv,\quad\quad  abuw$;$\quad\quad\quad\quad\quad  au(cv-bw) $ 

$AF : \quad acv^2, \quad\quad  abvw$;$\quad\quad\quad\quad\quad  av(cv-bw) $

$BA : \quad a^2vw, \quad\quad  acuv$;$\quad\quad\quad\quad\quad  av(aw-cu) $ 

$CA : \quad abuv, \quad\quad  abuv$;$\quad\quad\quad\quad\quad\quad\quad  0 $

$DA : \quad abvw, \quad\quad  bcuv$;$\quad\quad\quad\quad\quad  bv(aw-cu) $ 

$EA : \quad acuv, \quad\,  acuv$;$\quad\quad\quad\quad\quad\quad\quad  0 $

$FA : \quad acv^2, \quad\quad  bcuv$;$\quad\quad\quad\quad\quad  cv(av-bu) $

$BC : \quad abuw, \quad\quad  acuv$;$\quad\quad\quad\quad\quad  au(bw-cv) $ 

$BD : \quad abw^2, \quad\quad  acvw$;$\quad\quad\quad\quad\quad  aw(bw-cv) $

$BE : \quad acuw, \quad\quad  acuw$;$\quad\quad\quad\quad\quad\quad\quad  0 $

$BF :\quad acvw, \quad\quad  acvw$;$\quad\quad\quad\quad\quad\quad\quad  0 $

$CB : \quad abuw, \quad\quad  abuw$;$\quad\quad\quad\quad\quad\quad\quad  0 $ 

$DB : \quad abw^2, \quad\quad  bcuw$;$\quad\quad\quad\quad\quad  bw(aw-cu) $ 

$EB : \quad acuw, \quad\quad  acuw$;$\quad\quad\quad\quad\quad\quad\quad  0 $

$FB : \quad acvw, \quad\quad  bcuw$;$\quad\quad\quad\quad\quad  cw(av-bu) $

$CD : \quad b^2uw, \quad\quad  abvw$;$\quad\quad\quad\quad\quad  bw(bu-av) $

$CE : \quad bcu^2, \quad\quad  abuw,$;$\quad\quad\quad\quad\quad  bu(cu-aw) $

$CF : \quad bcuv, \quad\quad  abvw$;$\quad\quad\quad\quad\quad  bv(cu-aw) $

$DC : \quad b^2uw, \quad\quad  bcuv$;$\quad\quad\quad\quad\quad  bu(bw-cv) $

$EC : \quad bcu^2, \quad\quad  acuv$;$\quad\quad\quad\quad\quad  cu(bu-av) $

$FC : \quad bcuv, \quad\quad  bcuv$;$\quad\quad\quad\quad\quad\quad\quad  0 $

$DE : \quad bcuw, \quad\quad  bcuw$;$\quad\quad\quad\quad\quad\quad\quad  0 $

$DF : \quad bcvw, \quad\quad  bcvw$;$\quad\quad\quad\quad\quad\quad\quad  0 $

$ED : \quad bcuw, \quad\quad  acvw$;$\quad\quad\quad\quad\quad  cw(bu-av) $

$FD : \quad bcvw, \quad\quad  bcvw$;$\quad\quad\quad\quad\quad\quad\quad  0 $

$EF : \quad c^2uv, \quad\quad acvw$;$\quad\quad\quad\quad\quad  cv(cu-aw) $

$FE : \quad c^2uv, \quad\quad  bcuw$;$\quad\quad\quad\quad\quad  cu(cv-bw) $
\vspace{3mm}

From the list of coefficients above, we introduce the following notation:
\begin{align*}
[(xy)^2-x^2y^2]_1&:=av(bu-av)ijA+aw(cu-aw)ikB+bu(av-bu)ijC\\
&\quad\quad +bw(cv-bw)jkD+cu(aw-cu)ikE+cv(bw-cv)jkF,
\end{align*}
\begin{align*}
[(xy)^2-x^2y^2]_2&:=au(bu-av)ieC+aw(bu-av)ieD+au(cu-aw)ieE\\
&\quad\quad+av(cu-aw)ieF+bv(av-bu)jfA+bw(av-bu)jfB\\
&\quad\quad +bu(cv-bw)jfE+bv(cv-bw)jfF+cv(aw-cu)kgA\\
&\quad\quad+cw(aw-cu)kgB+cu(bw-cv)kgC+cw(bw-cv)kgD,
\end{align*}
\begin{align*}
[(xy)^2-x^2y^2]_3&:=au(av-bu)iAe+au(aw-cu)iBe+bu(aw-cu)iDe\\
&\quad\quad +cu(av-bu)iFe+av(bw-cv)jBf+bv(bu-av)jCf\\
&\quad\quad +bv(bw-cv)jDf+cv(bu-av)jEf+aw(cv-bw)kAg\\
&\quad\quad +bw(cu-aw)kCg+cw(cu-aw)kEg+cw(cv-bw)kFg,
\end{align*}
\begin{align*}
[(xy)^2-x^2y^2]_4&:=av(av-bu)A^2+aw(aw-cu)B^2+bu(bu-av)C^2\\
&\quad+bw(bw-cv)D^2+cu(cu-aw)E^2+cv(cv-bw)F^2,
\end{align*}
and
\begin{align*}
&[(xy)^2-x^2y^2]_5:=\\
&\quad aw(av-bu)AB+au(cv-bw)AE+av(cv-bw)AF+av(aw-cu)BA\\
&\quad +bv(aw-cu)DA+cv(av-bu)FA+au(bw-cv)BC+aw(bw-cv)BD\\
&\quad+bw(aw-cu)DB+cw(av-bu)FB+bw(bu-av)CD+bu(cu-aw)CE\\
&\quad +bv(cu-aw)CF+bu(bw-cv)DC+cu(bu-av)EC+cw(bu-av)ED\\
&\quad+cv(cu-aw)EF+cu(cv-bw)FE. 
\end{align*}
By simple calculations,  we obtain
\begin{equation} 
(xy)^2-x^2y^2=\sum_{\ell=1}^5[(xy)^2-x^2y^2]_\ell.
\end{equation}

For clarity and readability, we introduce the following notation:
\vspace{2mm}

$\alpha:=aw(av-bu), \beta:=au(cv-bw), \gamma:=av(cv-bw), \delta:=av(aw-cu)$,
\vspace{2mm}

$\varepsilon:=bv(aw-cu), \zeta:=cv(av-bu), \eta:=aw(bw-cv), \theta:=bw(aw-cu)$,
\vspace{2mm}

$\iota:=cw(av-bu), \kappa:=bw(bu-av), \lambda:=bu(cu-aw), \mu:=bu(bw-cv)$,
\vspace{2mm}

$\nu:=cu(bu-av), \xi:=cv(cu-aw), \pi:=cu(cv-bw)$.
\vspace{2mm}

Then we have
\begin{align*}
[(xy)^2-x^2y^2]_2&:=au(bu-av)ieC-\alpha ieD+au(cu-aw)ieE-\delta ieF\\
&\quad\quad+bv(av-bu)jfA-\kappa jfB-\mu jfE+bv(cv-bw)jfF\\
&\quad\quad -\xi kgA+cw(aw-cu)kgB-\pi kgC+cw(bw-cv)kgD,
\end{align*}

\begin{align*}
[(xy)^2-x^2y^2]_3&:=au(av-bu)iAe+au(aw-cu)iBe-\lambda iDe-\nu iFe\\
&\quad\quad -\gamma jBf+bv(bu-av)jCf+bv(bw-cv)jDf-\zeta jEf\\
&\quad\quad -\eta kAg-\theta kCg+cw(cu-aw)kEg+cw(cv-bw)kFg,
\end{align*}
and
\begin{align*}
[(xy)^2-x^2y^2]_5&:=\alpha AB+\beta (AE-BC)+\gamma AF+\delta BA+\varepsilon (DA-CF)\\
&\quad +\zeta FA+\eta BD+\theta DB+\iota (FB-ED)+\kappa CD+\lambda CE\\
&\quad +\mu DC+\nu EC+\xi EF+\pi FE. 
\end{align*}

Here we first assume that $\mathcal A$ is endo-commutative.
\vspace{2mm}

(i) Let $a=1$ and $b=c=0$.  Then
\[
\left\{\begin{array}{@{\,}lll} 
\alpha=\delta=vw,\\
\beta=\gamma=\varepsilon=\zeta=\eta=\theta=\iota=\kappa=\lambda=\mu=\nu=\xi=\pi=0,
\end{array} \right.
\]
and hence we have the following five equalities:
\vspace{2mm}

$[(xy)^2-x^2y^2]_1=-v^2ijA-w^2ikB$,

$[(xy)^2-x^2y^2]_2=-uvieC-vwie(D+F)-uwieE$,

$[(xy)^2-x^2y^2]_3=uviAe+uwiBe$,

$[(xy)^2-x^2y^2]_4=v^2A^2+w^2B^2$,

$[(xy)^2-x^2y^2]_5=vw(AB+BA)$.
\vspace{2mm}

\noindent
By (3), we have
\[
(xy)^2-x^2y^2=v^2(A^2-ijA)+w^2(B^2-ikB)+uvi(Ae-eC)
\]
\[
\quad\quad\quad\quad\quad\quad\quad\quad +vw\{AB+BA-ie(D+F)\}+uwi(Be-eE).
\]
Since $\mathcal A$ is endo-commutative by the assumption, it follows that
\begin{equation} 
\left\{\begin{array}{@{\,}lll} 
A^2=ijA,\\
B^2=ikB,\\
AB+BA=ie(D+F),\\
iAe=ieC,\\
iBe=ieE.
\end{array} \right.
\end{equation}

(ii) Let $b=1$ and $ a=c=0$.  Then
\[
\left\{\begin{array}{@{\,}lll} 
\kappa=\mu=uw,\\
\alpha=\beta=\gamma=\delta=\varepsilon=\zeta=\eta=\theta=\iota=\lambda=\nu=\xi=\pi=0, 
\end{array} \right.
\]
and hence we have the following five equalities:
\vspace{2mm}

$[(xy)^2-x^2y^2]_1=-u^2ijC-w^2jkD$,

$[(xy)^2-x^2y^2]_2=-uvjfA-uwjf(B+E)-vwjfF$,

$[(xy)^2-x^2y^2]_3=uvjCf+vwjDf$,

$[(xy)^2-x^2y^2]_4=u^2C^2+w^2D^2$,

$[(xy)^2-x^2y^2]_5=uw(CD+DC)$.
\vspace{2mm}

\noindent
By (3), we have
\[
(xy)^2-x^2y^2=u^2(C^2-ijC)+w^2(D^2-jkD)+uvj(Cf-fA)
\]
\[
\quad\quad\quad\quad\quad\quad\quad\quad +uw\{CD+DC-jf(B+E)\}+vwj(Df-fF),
\]
and hence
\begin{equation} 
\left\{\begin{array}{@{\,}lll} 
C^2=ijC,\\
D^2=jkD,\\
CD+DC=jf(B+E),\\
jCf=jfA,\\
jDf=jfF.
\end{array} \right.
\end{equation}

(iii) Let $c=1$ and $a=b=0$.  Then
\[
\left\{\begin{array}{@{\,}lll} 
\xi=\pi=uv,\\
\alpha=\beta=\gamma=\delta=\varepsilon=\zeta=\eta=\theta=\iota=\kappa=\lambda=\mu=\nu=0,
\end{array} \right.
\]
and hence we have the following five equalities:
\vspace{2mm}

$[(xy)^2-x^2y^2]_1=-u^2ikE-v^2jkF$,

$[(xy)^2-x^2y^2]_2=-uvkg(A+C)-uwkgB-vwkgD$,

$[(xy)^2-x^2y^2]_3=uwkEg+vwkFg$,

$[(xy)^2-x^2y^2]_4=u^2E^2+v^2F^2$,

$[(xy)^2-x^2y^2]_5=uv(EF+FE)$.
\vspace{2mm}

\noindent
By (3), we have
\[
(xy)^2-x^2y^2=u^2(E^2-ikE)+v^2(F^2-jkF)+uv\{EF+FE-kg(A+C)\}
\]
\[
\quad\quad\quad\quad\quad\quad\quad\quad +uwk(Eg-gB)+vwk(Fg-gD),
\]
and hence
\begin{equation} 
\left\{\begin{array}{@{\,}lll} 
E^2=ikE,\\
F^2=jkF,\\
EF+FE=kg(A+C),\\
kEg=kgB,\\
kFg=kgD.
\end{array} \right.
\end{equation}

(iv) Let $u=1$ and $v=w=0$.  Then
\[
\left\{\begin{array}{@{\,}lll} 
\lambda=\nu=bc,\\
\alpha=\beta=\gamma=\delta=\varepsilon=\zeta=\eta=\theta=\iota=\kappa=\mu=\xi=\pi=0,
\end{array} \right.
\]
and hence we have the following five equalities:
\vspace{2mm}

$[(xy)^2-x^2y^2]_1=-b^2ijC-c^2ikE$,

$[(xy)^2-x^2y^2]_2=abieC+acieE$,

$[(xy)^2-x^2y^2]_3=-abiAe-aciBe-bciDe-bciFe$,

$[(xy)^2-x^2y^2]_4=b^2C^2+c^2E^2$,

$[(xy)^2-x^2y^2]_5=bc(CE+EC)$.
\vspace{2mm}

\noindent
By (3), we have
\[
(xy)^2-x^2y^2=b^2(C^2-ijC)+c^2(E^2-ikE)+abi(eC-Ae)
\]
\[
\quad\quad\quad\quad\quad\quad\quad\quad +aci(eE-Be)+bc\{CE+EC-i(D+F)e\},
\]
and hence
\begin{equation} 
\left\{\begin{array}{@{\,}lll} 
C^2=ijC,\\
E^2=ikE,\\
CE+EC=i(D+F)e,\\
ieC=iAe,\\
ieE=iBe.
\end{array} \right.
\end{equation}

(v) Let $v=1$ and $u=w=0$.  Then
\[
\left\{\begin{array}{@{\,}lll} 
\gamma=\zeta=ac,\\
\alpha=\beta=\delta=\varepsilon=\eta=\theta=\iota=\kappa=\lambda=\mu=\nu=\xi=\pi=0,
\end{array} \right.
\]
and hence we have the following five equalities:
\vspace{2mm}

$[(xy)^2-x^2y^2]_1=-a^2ijA-c^2jkF$,

$[(xy)^2-x^2y^2]_2=abjfA+bcjfF$,

$[(xy)^2-x^2y^2]_3=-acjBf-abjCf-bcjDf-acjEf$,

$[(xy)^2-x^2y^2]_4=a^2A^2+c^2F^2$,

$[(xy)^2-x^2y^2]_5=ac(AF+FA)$.
\vspace{2mm}

\noindent
By (3), we have
\[
(xy)^2-x^2y^2=a^2(A^2-ijA)+c^2(F^2-jkF)+abj(fA-Cf)
\]
\[
\quad\quad\quad\quad\quad\quad\quad\quad +bcj(fF-Df)+ac\{AF+FA-j(B+E)f\},
\]
and hence
\begin{equation} 
\left\{\begin{array}{@{\,}lll} 
A^2=ijA,\\
F^2=jkF,\\
AF+FA=j(B+E)f,\\
jfA=jCf,\\
jfF=jDf.
\end{array} \right.
\end{equation}

(vi) Let $w=1$ and $u=v=0$.  Then
\[
\left\{\begin{array}{@{\,}lll} 
\eta=\theta=ab,\\
\alpha=\beta=\gamma=\delta=\varepsilon=\zeta=\iota=\kappa=\lambda=\mu=\nu=\xi=\pi=0,
\end{array} \right.
\]
and hence we have the following five equalities:
\vspace{2mm}

$[(xy)^2-x^2y^2]_1=-a^2ikB-b^2jkD$,

$[(xy)^2-x^2y^2]_2=ackgB+bckgD$,

$[(xy)^2-x^2y^2]_3=-abkAg-abkCg-ackEg-bckFg$,

$[(xy)^2-x^2y^2]_4=a^2B^2+b^2D^2$,

$[(xy)^2-x^2y^2]_5=ab(BD+DB)$.
\vspace{2mm}

\noindent
By (3), we have
\[
(xy)^2-x^2y^2=a^2(B^2-ikB)+b^2(D^2-jkD)+ack(gB-Eg)
\]
\[
\quad\quad\quad\quad\quad\quad\quad\quad +ab\{BD+DB-k(A+C)g\}+bck(gD-Fg),
\]
and hence
\begin{equation} 
\left\{\begin{array}{@{\,}lll} 
B^2=ikB,\\
D^2=jkD,\\
BD+DB=k(A+C)g,\\
kgB=kEg,\\
kgD=kFg.
\end{array} \right.
\end{equation}
Therefore, by combining equations (4) through (9), we obtain the following condition:
\begin{equation} 
\left\{\begin{array}{@{\,}lll} 
A^2=ijA, \cdots (10-1)\\
B^2=ikB,  \cdots (10-2)\\
C^2=ijC,  \cdots (10-3)\\
D^2=jkD,  \cdots (10-4)\\
E^2=ikE,  \cdots (10-5)\\
F^2=jkF,  \cdots (10-6)\\
AB+BA=ie(D+F),  \cdots (10-7)\\
CE+EC=i(D+F)e,  \cdots (10-8)\\
CD+DC=jf(B+E),  \cdots (10-9)\\
AF+FA=j(B+E)f,  \cdots (10-10)\\
EF+FE=kg(A+C),  \cdots (10-11)\\
BD+DB=k(A+C)g,  \cdots (10-12)\\
iAe=ieC,  \cdots (10-13)\\
iBe=ieE,  \cdots (10-14)\\
jCf=jfA,   \cdots (10-15)\\
jDf=jfF,  \cdots (10-16)\\
kEg=kgB,  \cdots (10-17)\\
kFg=kgD.  \cdots (10-18)
\end{array} \right.
\end{equation}
Consequently, we see that if $\mathcal A$ is endo-commutative, then (10) holds.
\vspace{2mm}

Conversely, we assume that the condition (10) holds.  Recall that
\begin{align*}
[(xy)^2-x^2y^2]_1&:=av(bu-av)ijA+aw(cu-aw)ikB+bu(av-bu)ijC\\
&\quad\quad +bw(cv-bw)jkD+cu(aw-cu)ikE+cv(bw-cv)jkF.
\end{align*}
By (10-1) to (10-6), we have
\begin{align*}
[(xy)^2-x^2y^2]_4&:=av(av-bu)A^2+aw(aw-cu)B^2+bu(bu-av)C^2\\
&\quad +bw(bw-cv)D^2+cu(cu-aw)E^2+cv(cv-bw)F^2\\
&\, \, =av(av-bu)ijA+aw(aw-cu)ikB+bu(bu-av)ijC\\
&\quad +bw(bw-cv)jkD+cu(cu-aw)ikE+cv(cv-bw)jkF.
\end{align*}
Therefore, we obtain
\begin{equation} 
[(xy)^2-x^2y^2]_1+[(xy)^2-x^2y^2]_4=0.
\end{equation}
By (10-13) to (10-18), we have
\begin{align*}
[(xy)^2-x^2y^2]_2&:=au(bu-av)ieC+aw(bu-av)ieD+au(cu-aw)ieE\\
&\quad\quad +av(cu-aw)ieF+bv(av-bu)jfA+bw(av-bu)jfB\\
&\quad\quad +bu(cv-bw)jfE+bv(cv-bw)jfF+cv(aw-cu)kgA\\
&\quad\quad+cw(aw-cu)kgB+cu(bw-cv)kgC+cw(bw-cv)kgD\\
&\, \, =au(bu-av)ieC-\alpha ieD+au(cu-aw)ieE-\delta ieF\\
&\quad\quad +bv(av-bu)jfA-\kappa jfB-\mu jfE+bv(cv-bw)jfF\\
&\quad\quad -\xi kgA+cw(aw-cu)kgB-\pi kgC+cw(bw-cv)kgD\\
&\, \, =au(bu-av)iAe+au(cu-aw)iBe+bv(av-bu)jCf\\
&\quad\quad +bv(cv-bw)jDf+cw(aw-cu)kEg+cw(bw-cv)kFg\\
&\quad\quad  -\alpha ieD-\delta ieF-\kappa jfB-\mu jfE-\xi kgA-\pi kgC.
\end{align*}
Recall that
\begin{align*}
[(xy)^2-x^2y^2]_3&:=au(av-bu)iAe+au(aw-cu)iBe-\lambda iDe-\nu iFe\\
&\quad\quad -\gamma jBf+bv(bu-av)jCf+bv(bw-cv)jDf-\zeta jEf\\
&\quad\quad -\eta kAg-\theta kCg+cw(cu-aw)kEg+cw(cv-bw)kFg,
\end{align*}
Therefore, we obtain
\begin{equation} 
[(xy)^2-x^2y^2]_2+[(xy)^2-x^2y^2]_3=-\alpha ieD-\delta ieF-\kappa jfB-\mu jfE-\xi kgA-\pi kgC
\end{equation}
\[
\quad\quad\quad\quad\quad\quad\quad\quad\quad\quad\quad\quad\quad\quad\quad\quad\quad -\lambda iDe-\nu iFe-\gamma jBf-\zeta jEf-\eta kAg-\theta kCg.
\]
Recall that
\begin{equation} 
[(xy)^2-x^2y^2]_5=\alpha AB+\beta(AE-BC)+\gamma AF+\delta BA+\varepsilon(DA-CF)+\zeta FA
\end{equation}
\[
\quad\quad\quad\quad\quad\quad\quad\quad +\, \eta BD+\theta DB+\iota(FB-ED)+\kappa CD+\lambda CE+\mu DC
\]
$\quad\quad\quad\quad\quad\quad\quad\quad\quad\quad\, \, +\, \nu EC+\xi EF+\pi FE$.
\vspace{1mm}

\noindent
By (3),  (11), (12) and (13), we have
\vspace{1mm}

\begin{equation} 
(xy)^2-x^2y^2=\alpha(AB-ieD)+\beta(AE-BC)+\gamma(AF-jBf)+\delta(BA-ieF)+\varepsilon(DA-CF)
\end{equation}
$\quad\quad\quad\quad\quad\quad\quad\quad\quad+\, \zeta(FA-jEf)+\eta(BD-kAg)+\theta(DB-kCg)+\iota(FB-ED)$

$\quad\quad\quad\quad\quad\quad\quad\quad +\, \kappa(CD-jfB)+\lambda(CE-iDe)+\mu(DC-jfE)$

$\quad\quad\quad\quad\quad\quad\quad\quad+\, \nu(EC-iFe)+\xi(EF-kgA)+\pi(FE-kgC)$.
\vspace{2mm}

Now, we can easily verify that
\begin{equation}
\left\{\begin{array}{@{\,}lll} 
\alpha-\delta=\beta, \\
\gamma-\zeta=-\varepsilon, \\
\eta-\theta=-\iota, \\
\kappa-\mu=-\varepsilon, \\
\lambda-\nu=\beta, \\
\xi-\pi=-\iota.
\end{array} \right.
\end{equation}
By (15) and (10-7), we have
\[
A_1:=\alpha(AB-ieD)+\beta(AE-BC)+\gamma(AF-jBf)+\delta(BA-ieF)+\varepsilon(DA-CF)
\]
\[
=\alpha(AB-ieD)+\beta(AE-BC)+\gamma(AF-jBf)+(\alpha-\beta)(BA-ieF)+\varepsilon(DA-CF)
\]
\[
=\alpha AB-\alpha ieD+\beta(AE-BC)+\gamma(AF-jBf)+\alpha(BA-ieF)-\beta(BA-ieF)+\varepsilon(DA-CF)
\]
\[
=\alpha(AB+BA)-\alpha ieD+\beta(AE-BC)+\gamma(AF-jBf)-\alpha ieF-\beta(BA-ieF)+\varepsilon(DA-CF)
\]
\[
=\alpha ie(D+F)-\alpha ieD+\beta(AE-BC)+\gamma(AF-jBf)-\alpha ieF-\beta(BA-ieF)+\varepsilon(DA-CF)
\]
\[
=\beta(AE-BC)+\gamma(AF-jBf)-\beta(BA-ieF)+\varepsilon(DA-CF)
\]

\noindent
By (15) and (10-12), we have
\vspace{2mm}

$A_2:=\zeta(FA-jEf)+\eta(BD-kAg)+\theta(DB-kCg)+\iota(FB-ED)$

\quad\,\,  $=\zeta(FA-jEf)+\eta(BD-kAg)+(\eta+\iota)(DB-kCg)+\iota(FB-ED)$


\quad\,\,  $=\zeta(FA-jEf)+\eta BD-\eta kAg+\eta DB+\iota DB-(\eta+\iota)kCg+\iota(FB-ED)$

\quad\,\,  $=\zeta(FA-jEf)+\eta(BD+DB)-\eta kAg-(\eta+\iota) kCg+\iota(DB+FB-ED)$

\quad\,\,  $=\zeta(FA-jEf)+\eta k(A+C)g-\eta kAg-\theta kCg+\iota(DB+FB-ED)$

\quad\,\,  $=\zeta(FA-jEf)+\eta kCg-\theta kCg+\iota(DB+FB-ED)$

\quad\,\,  $=\zeta(FA-jEf)+(\eta-\theta)kCg+\iota(DB+FB-ED)$

\quad\,\,  $=(\gamma+\varepsilon)(FA-jEf)-\iota kCg+\iota(DB+FB-ED)$.
\vspace{2mm}

\noindent
By (15) and (10-9), we have
\vspace{2mm}

$A_3:=\kappa(CD-jfB)+\lambda(CE-iDe)+\mu(DC-jfE)$

\quad\, \, $=\kappa(CD-jfB)+\lambda(CE-iDe)+(\kappa+\varepsilon)(DC-jfE)$

\quad\, \, $=\kappa CD-\kappa jfB+\lambda(CE-iDe)+\kappa DC-\kappa jfE+\varepsilon(DC-jfE)$

\quad\, \, $=\kappa(CD+DC)-\kappa jfB+\lambda(CE-iDe)-\kappa jfE+\varepsilon(DC-jfE)$

\quad\, \, $=\kappa jf(B+E)-\kappa jfB+\lambda(CE-iDe)-\kappa jfE+\varepsilon(DC-jfE)$

\quad\, \, $=\lambda(CE-iDe)+\varepsilon(DC-jfE)$
\vspace{2mm}

\noindent
By (15) and (10-11), we have
\vspace{2mm}

$A_4:=\nu(EC-iFe)+\xi(EF-kgA)+\pi(FE-kgC)$

\quad\, \, $=\nu(EC-iFe)+\xi(EF-kgA)+(\xi+\iota)(FE-kgC)$

\quad\, \, $=\nu(EC-iFe)+\xi EF-\xi kgA+\xi FE-\xi kgC+\iota(FE-kgC)$

\quad\, \, $=\nu(EC-iFe)+\xi(EF+FE)-\xi kgA-\xi kgC+\iota(FE-kgC)$

\quad\, \, $=\nu(EC-iFe)+\xi kg(A+C)-\xi kgA-\xi kgC+\iota(FE-kgC)$

\quad\, \, $=(\lambda-\beta)(EC-iFe)+\iota(FE-kgC)$.
\vspace{2mm}

\noindent
Therefore, by (10-10) and (10-8), we have
\begin{align*}
A_1+&A_2+A_3+A_4\\
&=\beta(AE-BC)+\gamma(AF-jBf)-\beta(BA-ieF)+\varepsilon(DA-CF)\\
&\quad +(\gamma+\varepsilon)(FA-jEf)-\iota kCg+\iota(DB+FB-ED)\\
&\quad +\lambda(CE-iDe)+\varepsilon(DC-jfE)+(\lambda-\beta)(EC-iFe)+\iota(FE-kgC)\\
&=\beta(AE-BC)+\gamma(AF+FA)-\gamma jBf-\beta(BA-ieF)+\varepsilon(DA-CF)\\
&\quad +\varepsilon FA-(\gamma+\varepsilon)jEf-\iota kCg+\iota(DB+FB-ED)\\
&\quad +\lambda(CE+EC)-\lambda iDe+\varepsilon(DC-jfE)-\beta EC-(\lambda-\beta)iFe+\iota(FE-kgC)\\
&=\beta(AE-BC)+\gamma j(B+E)f-\gamma jBf-\beta(BA-ieF)+\varepsilon(DA-CF)\\
&\quad +\varepsilon FA-(\gamma+\varepsilon)jEf-\iota kCg+\iota(DB+FB-ED)\\
&\quad +\lambda i(D+F)e-\lambda iDe+\varepsilon(DC-jfE)-\beta EC-(\lambda-\beta)iFe+\iota(FE-kgC)\\
&=\beta(AE-BC)-\beta(BA-ieF)+\varepsilon(DA-CF)\\
&\quad +\varepsilon FA-\varepsilon jEf-\iota kCg+\iota(DB+FB-ED)\\
&\quad +\varepsilon(DC-jfE)-\beta EC+\beta iFe+\iota(FE-kgC)\\
&=\beta(AE-BC-BA+ieF-EC+iFe)\\
&\quad +\varepsilon(DA-CF+FA-jEf+DC-jfE)\\
&\quad +\iota(-kCg+DB+FB-ED+FE-kgC).
\end{align*}
However, since $(xy)^2-x^2y^2=A_1+A_2+A_3+A_4$ from (14), we have
\begin{equation}
(xy)^2-x^2y^2=\beta(AE-BC-BA+ieF-EC+iFe)
\end{equation} 
\[
\quad\quad\quad\quad\quad\quad\quad\quad\quad +\varepsilon(DA-CF+FA-jEf+DC-jfE)
\]
\[
\quad\quad\quad\quad\quad\quad\quad\quad\quad\quad\quad\quad +\iota(-kCg+DB+FB-ED+FE-kgC).
\]
\vspace{2mm}

In light of equation (16), consider the following condition:
\begin{equation} 
\left\{\begin{array}{@{\,}lll} 
BC+BA+EC-AE=i(eF+Fe), \\
DA+FA+DC-CF=j(Ef+fE), \\
DB+FB+FE-ED=k(Cg+gC). 
\end{array} \right.
\end{equation}
If (17) holds, we see from (16) that $(xy)^2=x^2y^2$ for all $x, y\in \mathcal A$, that is, $\mathcal A$ is endo-commutative.  

Consequently, we see that if both (10) and (17) hold, then $\mathcal A$ is endo-commutative.  
\vspace{2mm}

\begin{thm} 
Let $\mathcal A$ be a 3-dimensional curled algebra of type $(i, j, k)$ over $K$ with linear basis $\{e, f, g\}$.  Then $\mathcal A$ is endo-commutative iff (10) and (17) hold.
\end{thm}

\begin{proof}
%
%
%
%
%
  By the argument immediately preceding the theorem, we have shown that if conditions (10) and (17) hold, then $\mathcal A$ is endo-commutative.
  Conversely, we have already proved that if $\mathcal A$ is endo-commutative, then (10) holds. Moreover, we have demonstrated that
  (10) implies (16). Namely, if $\mathcal A$ is endo-commutative, then \[\beta(AE-BC-BA+ieF-EC+iFe)+\varepsilon(DA-CF+FA-jEf+DC-jfE)\]
  \[\quad\quad\quad\quad\quad\quad\quad+\iota(-kCg+DB+FB-ED+FE-kgC)=0\]
  for all $a, b, c, u, v, w\in K$. We now consider the following three cases:

  (i) Let $b=w=0, a=u=c=v=1$. Then, $\beta=1, \varepsilon=\iota=0$, and hence \[AE-BC-BA+ieF-EC+iFe=0.\]

  (ii) Let $c=u=0, b=v=a=w=1$. Then, $\varepsilon=1, \beta=\iota=0$, and hence \[DA-CF+FA-jEf+DC-jfE=0.\]

  (iii) Let $b=u=0, c=w=a=v=1$. Then, $\iota=1, \beta=\varepsilon=0$, and hence \[-kCg+DB+FB-ED+FE-kgC=0.\]

  Therefore, condition (17) also holds.
  
\end{proof}

Remark.  Let $\mathcal A$ be a 3-dimensional curled algebra of type $(i, j, k)$ over $K$ with linear basis $\{e, f, g\}$.  Then we see easily that $\mathcal A$ is zeropotent iff
\begin{equation} 
\left\{\begin{array}{@{\,}lll} 
i=j=k=0, \\
A+C=0, \\
B+E=0, \\
D+F=0
\end{array} \right.
\end{equation}
holds. The reader may find it helpful to compare conditions (10) and (17) with condition (18).

\end{document}